\documentclass[11pt,twoside]{article}
\usepackage{amssymb, amsmath, epsfig, amsfonts, pst-plot, pst-node, pst-text, setspace, sectsty, natbib}

\newtheorem{theorem}{Theorem}
\newtheorem{lemma}{Lemma}

\newtheorem{corollary}{Corollary}

\newcounter{cremark}
\newenvironment{remark}[0]{\refstepcounter{cremark}\begin{trivlist}
\item[\hskip \labelsep {\emph{Remark
\thecremark.}}]}{\end{trivlist}}

\newenvironment{proof}[0]{\begin{trivlist}
\item[\hskip \labelsep {\emph{Proof.}}]}{\end{trivlist}}

\newenvironment{proof_of}[1][]{\begin{trivlist}
\item[\hskip \labelsep {\bfseries Proof of #1:}]}{\end{trivlist}}

\newcounter{cexample}
\newenvironment{example}[0]{\refstepcounter{cexample}\begin{trivlist}
\item[\hskip \labelsep {\emph{Example
\thecexample.}}]}{\end{trivlist}}

\newcommand{\qed}{\nobreak \ifvmode \relax \else
      \ifdim\lastskip<1.5em \hskip-\lastskip
      \hskip1.5em plus0em minus0.5em \fi \nobreak
      $\Box$\fi}

\newcommand{\ind}[1]{1\{ {#1} \} }

\newcommand{\R}[0]{\mathbb{R}}

\renewcommand{\H}[0]{\mathbb{H}}
\newcommand{\N}[0]{\mathbb{N}}

\newcommand{\V}[0]{\mathbb{V}}

\renewcommand{\P}[0]{\mathbb{P}}
\newcommand{\Rplus}[0]{\mathbb{R}_+}

\newcommand{\st}[0]{:}

\def\convas{\rightarrow_{\text{a.s.}}}

\def\tilde{\widetilde}
\def\hat{\widehat}

 \allsectionsfont{\normalsize}

\makeatletter
\def\@seccntformat#1{\csname the#1\endcsname.\quad}
\makeatother

\usepackage{geometry}%
\begin{document}

\title{\normalfont\Large\mdseries\bf{Inconsistency of the MLE
for the joint distribution of interval censored survival times
 and continuous marks}}
\author{\sc{M.H. Maathuis and J.A. Wellner}\\
\emph{Department of Statistics, University of Washington}\\
\date{}
} \maketitle

\vspace{-.5cm}

\noindent ABSTRACT. This paper considers the nonparametric maximum
likelihood estimator (MLE) for the joint distribution function of an
interval censored survival time and a continuous mark variable. We
provide a new explicit formula for the MLE in this problem. We use
this formula and the mark specific cumulative hazard function of
\citet*{HuangLouis98} to obtain the almost sure limit of the MLE. This result
leads to necessary and sufficient conditions for consistency of the MLE
which imply that the MLE is inconsistent in general. We show that the inconsistency
can be repaired by discretizing the marks. Our theoretical results are supported by
simulations.  \\


\section{Introduction}\label{sec: introduction}

Suppose that $X$ is a survival time and $Y$ is a continuous mark
variable which may be correlated with $X$. \citet{HuangLouis98}
considered nonparametric estimation of the joint distribution of $X$
and $Y$ when $X$ is subject to (random) right-censoring and the mark
variable $Y$ is observed if and only if $X$ is uncensored.  In many
cases of interest, however, we can only observe an interval censored
version of the random variable $X$. For example, \citet*[henceforth
HMG]{HudgensMaathuisGilbert07} analyzed an HIV vaccine trial in
which $X$ is the time of HIV infection and $Y$ is a measure of the
genetic distance between the infecting HIV virus and the virus in
the vaccine. The participants of this trial were tested for HIV at
several follow-up times. As a result, $X$ was interval censored,
that is, only known to be in a time interval determined by the
follow-up times. Moreover, since the viral distance $Y$ could only
be determined for HIV positive individuals, $Y$ was missing for all
individuals who were HIV negative at their last follow-up visit.

Motivated by this example we consider the following model, that we
refer to as the ``interval censored  continuous mark model". Let
$X>0$ be a survival time and let $Y \in \R$ be a continuous mark
variable. For a fixed integer $k\ge 1$, suppose that ${\bf T} = (T_1
, \ldots , T_k)$ is a vector of observation times with distribution
$G$. We assume that $0 < T_1 < \cdots < T_k $ and that ${\bf T}$ is
independent of $(X,Y)$. We cannot observe $(X,Y)$ directly. Instead,
our observed data are $W= ( {\bf T}, {\bf \Delta}, Z)$, where
\vspace{-.1cm}
\begin{align*}
  {\bf \Delta} = (\Delta_1 , \ldots , \Delta_{k+1} )
  \quad \text{with} \quad \Delta_j \equiv 1 \{ T_{j-1} < X \le T_j \}, \ \ j = 1, \ldots , k+1,
\end{align*}
(with the convention that $T_0 \equiv 0$ and $T_{k+1} \equiv \infty$), and
\vspace{-.1cm}
\begin{align*}
  Z = \Delta_+ Y  \quad \text{with} \quad \Delta_+ \equiv \sum_{j=1}^k \Delta_j.
\end{align*}
Note that the vectors ${\bf T}$ and ${\bf \Delta}$ determine a time interval $(T_{j-1},T_j]$, $j=1,\dots,k+1$,
that is known to contain the survival time $X$. The variable $Z$ reflects that the
mark variable $Y$ is observed if and only if the survival endpoint is reached
before the last observation time, i.e., if and only if $X\le T_k$.

Our censoring model for $X$ is called ``interval censoring case
$k$'', since each individual in the study has exactly $k$
observation times $T_1,\dots,T_k$
(see \citet*{GroeneboomWellner92} for case 1 and case 2
interval censoring, and \citet*{Wellner95} for case $k$ interval
censoring). Interval censoring case 1 is also referred to
as ``current status censoring", since we only observe the ``current
status" of an individual at a single observation time. A model
which allows the number of observation times to be random,
and hence to vary across individuals in the study, is called ``mixed
case interval censoring'' (see e.g. \citet*{SchickYu00},
\citet*{VanderVaartWellner00:GC-preservationTheorems}, and
\citet*[page 12]{Sun06}).

Our goal here is to study the nonparametric maximum likelihood estimator (MLE) of the joint distribution
$F_0$ of $(X,Y)$ when the observations consist of $W_1 , \ldots , W_n$ i.i.d. as $W$.
In particular we focus on consistency issues, and we show, in fact, that the MLE is inconsistent
in general.

There are several known examples of inconsistency of the nonparametric maximum likelihood estimator.
\citet[pages 255 -- 258]{BarlowEtAl72}
showed that the MLE $\widehat{F}_n$  for the class of
star-shaped distributions  (distributions on $[0,b)$ with $F(0) = 0$ and $F(x)/x$ non-decreasing)
is inconsistent, by showing that for sampling from the uniform distribution on $[0,1]$ the
MLE $\widehat{F}_n (x) \rightarrow_{a.s.} x^2$.
For distributions $F$ with increasing failure rate average (IFRA), \citet{BoylesMarshallProschan85} showed that the
MLE is inconsistent, and they identified the limit explicitly for sampling from a general continuous
distribution function $F$.
In the context of bivariate right-censored data, inconsistency of the nonparametric MLE
for continuous bivariate distributions was pointed out by \citet{TsaiLeurgansCrowley86} and was also studied by
\citet*{VanderLaan96}.
For estimation of a distribution function on $\R$ based on left-truncated and
case 1 interval censored data,
\citet{PanChappell99} showed that the nonparametric MLE is inconsistent. Finally, \citet*[Section
6.2]{Maathuis03} showed inconsistency of the MLE of the bivariate distribution of $(X,Y)$ when
$X$ is subject to current status censoring and $Y$ is observed exactly.

There are many more examples of inconsistent maximum likelihood estimators in parametric problems: see, for example,
\citet{NeymanScott48},
\citet{Bahadur58}, \citet{Ferguson82}, \citet{GhoshYang95}, \citet{GuptaSzekelyZsigiri99}, and the
interesting review by \citet{LeCam90}.

To relate our inconsistency result to some of these earlier studies
of inconsistency of the MLE, note that observation of $W$ instead of
$(X,Y)$ can be regarded as observation of a (random) set $A$ known
to contain the unobservable $(X,Y)$. We call such a set an
\emph{observed set}. In our model the observed sets can take two
forms. When $\Delta_j = 1$ for some $j \le k$ (so $\Delta_+=1$),
then the observed set is a horizontal line segment: \vspace{-.15cm}
\begin{align}\label{eq: observed set Delta+ = 1}
   A = (T_{j-1},T_j] \times \{Z\},
\end{align}
while when $\Delta_{k+1} = 1$, or equivalently, when $\Delta_+=0$, the observed
set is a half plane:
\vspace{-.15cm}
\begin{align}\label{eq: observed set Delta+ = 0}
  \qquad A = (T_{k}, \infty) \times \R.
\end{align}
The line segments that arise when $\Delta_+=1$ are an indicator of
potential consistency problems for the MLE,
since such line segments also occurred in the inconsistent MLEs studied by
\citet*{VanderLaan96} and \citet*[Section
6.2]{Maathuis03}. This prompted us
to carefully study consistency of the MLE
for interval censored continuous mark data.

Our work is also related to the classical competing risks model, in
which one studies the failure time $X$ of a system that can fail
from a (finite) number of $J$ competing risks given by values of $Y
\in \{1,\dots,J\}$. The variable $Y$ in this model can only be
observed after the failure event happened, and is therefore a mark
variable. Thus, the classical competing risks model can be called a
``discrete mark model", and can be viewed as the discrete
counterpart of the continuous mark model. The competing risks model
has been studied under various censoring assumptions for $X$.
\citet*{Aalen76, Aalen78} and \citet*[\S 7.2, pages 163 -
178]{KalbfleischPrentice80} studied the MLE in this model when $X$
is subject to right censoring. The generalization to interval
censored survival data with competing risks was considered by
\citet*{HudgensSattenLongini01} and
\citet*{JewellVanderLaanHenneman03}. \citet*{JewellKalbfleisch04}
studied computational issues of the MLE for current status data with
competing risks, and \citet*{Maathuis06thesis},
\citet*{GroeneboomMaathuisWellner06a}, and
\citet*{GroeneboomMaathuisWellner06b} derived the asymptotic
properties of the MLE in this model.

In the current paper we focus on the interval censored continuous mark model.
In Section~\ref{sec: closed form} we derive a new formula for the MLE in this model, using
connections with univariate right censored data. In Section~\ref{sec:limit} we use
this new formula and the mark specific cumulative hazard function of
\citet*{HuangLouis98} to derive the almost sure limit of the MLE.
This result leads to necessary and sufficient conditions for consistency of the MLE which
force a relation between the unknown distribution $F_0$ and the observation time distribution $G$.
Since such a relation will typically not hold, it follows that the MLE
is inconsistent in general. In Section~\ref{sec: repaired mle} we
show that the inconsistency can be repaired by discretizing the
marks, an operation that transforms the data into interval censored competing risks data.
In Section~\ref{sec: examples} we support our theoretical results by simulations
of the MLE and the repaired MLE.
Section~\ref{sec: discussion} contains a discussion of some remaining issues.
Technical proofs are collected in the Appendix, Section~\ref{sec: appendix}.

\section{Explicit formula for the MLE}\label{sec: closed form}

HMG noted a close connection between the MLE for univariate right
censored data and the MLE for interval censored continuous mark
data. We use this connection in Section~\ref{subsec:closedformcontinuousmark}
to derive a new explicit formula for the MLE
for interval censored continuous mark data. But first, in Section \ref{subsec: intermezzo},
we review univariate right censored data in a way that shows the
similarity between the two models.

\subsection{Intermezzo: univariate right censored data}\label{subsec: intermezzo}

Suppose that we want to estimate the distribution $F_0$ of a survival time $X$,
and suppose that $X$ is subject
to right censoring. Thus, instead of $n$ i.i.d.\ copies of $X$, we
observe $n$ i.i.d.\ copies of $(\min(X,T), \ind{X\le T})$, where $T$
is a random censoring time with distribution $G$. We assume that $T$ is independent of $X$. It is
well-known that the MLE $\hat F_n$ of $F_0$ in this model is given by the Kaplan-Meier estimator.

We now review the Kaplan-Meier estimator in a way that allows us to easily make
a connection with interval censored continuous mark data. We first
introduce some notation and terminology. Define $U \equiv \min(X,T)$
and $\Delta \equiv \ind{X\le T}$, and let $(U_1,\Delta_1),\dots,(U_n,\Delta_n)$
denote $n$ i.i.d.\ copies of $(U,\Delta)$. Recalling the discussion of
observed sets in Section \ref{sec: introduction}, each observation $(U,\Delta)$ defines
an observed set $A$ that is known to contain $X$: $A = \{U\}$ if $\Delta=1$,
and $A = (U,\infty)$ if $\Delta=0$. Let $U_{(1)},\dots,U_{(n)}$ be the
order statistics of $U_1,\dots,U_n$, and let $\Delta_{(i)}$ and
$A_{(i)}$ be the corresponding values of $\Delta$ and $A$. We assume
that all $A_i$ with $\Delta_i=1$ are distinct, since this will be
the case for the continuous mark data. However, we allow ties in the
$T$'s and $U$'s provided that this assumption is not violated. We break
such ties in $U$ arbitrarily after ensuring that observations with
$\Delta=1$ are ordered before those with $\Delta=0$.

By assuming that $F$ has a density $f$ with respect to some
dominating measure $\mu$, the likelihood
(up to multiplicative terms depending only on $G$) is
$L_n(F) = \prod_{i=1}^n q(U_i,\Delta_i)$, where $ q(u,\delta) =
f(u)^{\delta}\left\{1-F(u)\right\}^{1-\delta}$. Since the first term of
$q$ is a density-type term, $L_n(F)$ can be made
arbitrarily large by letting $f$ peak at some value $U_i$ with
$\Delta_i=1$. This problem is usually solved by maximizing $L_n(F)$
over the class of distribution functions that have a density with
respect to counting measure on the observed failure times. We can
then write $L_n(F) = \prod_{i=1}^n P_F(A_i)$, where $P_F(A)$ is the
probability of $A$ under $F$.

It is well-known \citep{Peto73, Turnbull76} that the MLE in censored
data problems can only assign mass to a finite number of disjoint
regions, called \emph{maximal intersections} by \citet*{WongYu99}.
\citet{Maathuis05} introduced an efficient algorithm to compute the
maximal intersections for $d$-variate interval censored data. This
algorithm is based on a height map $h : \R^d \to \N$ of the observed
sets, where $h(x)$ is defined as the number of observed sets that
contain $x$. Maathuis showed that the maximal intersections
correspond exactly to the local maximum regions of the height map of
the observed sets. (If there are ties in the observed sets, then
these need to be resolved before applying the height map, see
\citet*{Maathuis05}.)

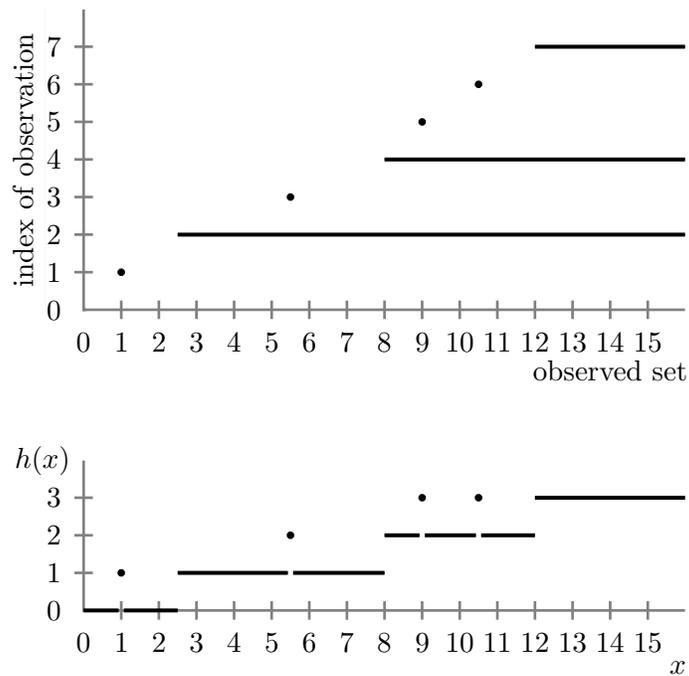
\begin{figure}
   \psset{unit=5mm}
   \psset{fillstyle=solid}
   \psset{fillcolor=black}
   \psset{linewidth=1.5pt}
   \begin{center}
         \begin{pspicture}(-1,-2)(16,16)
         \cnode[linecolor=white,fillcolor=white](-1,8){0}{aa}
         \cnode[linecolor=white,fillcolor=white](-1,15.99){0}{bb}
         \ncline[linewidth=0pt,linecolor=white]{aa}{bb}
         \naput{\rotateleft{index of observation}}
         \uput[d](14,7){observed set}
         \psaxes[linewidth=1pt,linecolor=gray](0,8)(15.99,15,99)
         \cnode(1,9){.1}{0}
         \psline(2.5,10)(15.99,10)
         \cnode(5.5,11){.1}{2}
         \psline(8,12)(15.99,12)
         \cnode(9,13){.1}{3}
         \cnode(10.5,14){.1}{4}
         \psline(12,15)(15.99,15)
         \uput[l](0,4){$h(x)$}
         \uput[d](15.8,-1){$x$}
         \psaxes[linecolor=gray,linewidth=1pt](0,0)(15.99,3.99)
         \psline(0,0)(.93,0)
         \cnode(1,1){.1}{A}
         \psline(1.07,0)(2.5,0)
         \psline(2.5,1)(5.43,1)
         \cnode(5.5,2){.1}{B}
         \psline(5.57,1)(8,1)
         \psline(8,2)(8.93,2)
         \cnode(9,3){.1}{C}
         \psline(9.07,2)(10.43,2)
         \cnode(10.5,3){.1}{D}
         \psline(10.57,2)(12,2)
         \psline(12,3)(15.99,3)
         \end{pspicture}
         \caption{Observed sets (upper panel) and the corresponding height map (lower panel)
         for univariate right censored data, based on the following 7 observations of $(U,\Delta)$:
         $(1,1)$, $(2.5,0)$, $(5.5,1)$, $(8,0)$, $(9,1)$, $(10.5,1)$ and $(12,0)$.
         Note that the maximal intersections are given by the local maximum regions of the height map:
         $\{1\}$, $\{5.5\}$, $\{9\}$, $\{10.5\}$ and $(12,\infty)$.}
         \label{fig: height map}
   \end{center}
\end{figure}

The height map $h: \R \mapsto \N$ for univariate right
censored data is illustrated in Figure \ref{fig: height map}. Note that
$h(x)$ simply represents the number of observed
sets $A_1,\dots,A_n$ that overlap at the point $x$. It is clear that all
sets $A_{(i)}$ with $i\in \mathcal I = \{i\in \{1,\dots,n\} \st
\Delta_{(i)} =1\}$, or in other words, all sets of the form $A_{(i)}=\{U_{(i)}\}$,
are local maxima of the height map. Hence, all such sets are
maximal intersections, and we denote these by $M_{(i)}$, $i\in \mathcal I$. This
notation may seem redundant since $M_{(i)} = A_{(i)}$, but it will
be useful in Section \ref{subsec:closedformcontinuousmark}. Furthermore, if and only if
$\Delta_{(n)}=0$, the height map has an
extra local maximum region $A_{(n)} = (U_{(n)},\infty)$, resulting in an extra
maximal intersection
$M_{(n+1)} = (U_{(n)},\infty)$. This situation is
illustrated in Figure \ref{fig: height map}.
Let $\overline{\mathcal I}$ be the collection of indices of all maximal
intersections. Thus, $\overline{\mathcal I} = \mathcal I$ if
$\Delta_{(n)}=1$ and $\overline{\mathcal I} = \mathcal I \cup
\{n+1\}$ if $\Delta_{(n)}=0$.

{\bf [Figure \ref{fig: height map} about here.]}

Let $p_i$ be the probability mass of
maximal intersection $M_{(i)}$, $i \in \overline{\mathcal I}$. We
can then write the likelihood in terms of the $p_i$'s:\vspace{-.1cm}
\begin{align} \label{eq: likelihood hood univariate right
censored}
   \prod_{i=1}^n P(A_i) = \prod_{i=1}^n  \left(\sum_{j\in
   \overline{\mathcal I}} p_j\ind{M_{(j)}\subseteq A_{(i)}}\right) =
   \prod_{i=1}^n p_i^{\Delta_{(i)}} \left( \sum_{j\ge i+1, j\in \overline{\mathcal
   I}} p_j \right)^{1-\Delta_{(i)}},
\end{align}
where the second equality follows from the fact that the data are ordered
with respect to the variable $U=\min(X,T)$.
The MLE $\hat p$ maximizes this expression under the constraints\vspace{-.1cm}
\begin{align}
   \sum_{i\in \overline{\mathcal I}} p_i = 1 \quad \text{and}
   \quad p_i \ge 0 \quad \text{for all}\; i\in \overline{\mathcal I}.\label{eq: constraints}
\end{align}
It is well-known that $\hat p$ is the Kaplan-Meier or product-limit
estimator, given by\vspace{-.1cm}
\begin{align*}
   \hat p_i & = \prod_{j=1}^{i-1} \left(1-\frac{\Delta_{(j)}}{n-j+1}\right)
   \frac{\Delta_{(i)}}{n-i+1},\qquad i \in \mathcal I,
\end{align*}
and $\hat p_{n+1}= 1-\sum_{i \in \mathcal I} \hat p_i$ if
$\Delta_{(n)}=0$ (see for example \citet*{ShorackWellner86}, Chapter
7, pages 332-333). Equivalently, we can write\vspace{-.1cm}
\begin{align*}
   \sum_{j\ge i, j\in \overline{\mathcal I} } \hat p_j = \prod_{j\le i-1} \left( 1-
   \frac{\Delta_{(j)}}{n-j+1}\right), \qquad
   i \in \overline{\mathcal I}.
\end{align*}
The vector $\hat p$ is uniquely determined. We obtain $\hat F_n(x)$
by summing all probability mass of $\hat p$ that falls in the interval $(0,x]$. It is
well-known that $\hat F_n(x)$ is non-unique for $x
> U_{(n)}$ if and only if $\Delta_{(n)}=0$. This is
caused by the fact that the MLE is indifferent to the distribution
of mass within a maximal intersection, called ``representational
non-uniqueness" by \citet*{GentlemanVandal02}. Since all maximal
intersections $\{M_{(i)}\st i\in \mathcal I\}$ are points, this
non-uniqueness occurs if and only if $M_{(n+1)} = (U_{(n)},\infty)$
exists, and this happens if and only if $\Delta_{(n)}=0$.

\subsection{Continuous mark data: Explicit formula for the MLE}\label{subsec:closedformcontinuousmark}

We now return to the interval censored continuous mark model given in Section \ref{sec: introduction}, and
introduce some additional notation.
Let $F_0(x,y) = P(X\le x, Y\le y)$ be the joint distribution of $(X,Y)$, and let
$F_{0X}(x) = F_{0} (x,\infty)=P(X\le x)$ and
$F_{0Y}(y)=F_{0} (\infty, y) = P(Y\le y)$ be the marginal distributions of $X$ and $Y$,
respectively. Recall that $G$ denotes the distribution of the observation times $\mathbf{T}$.
We use subscripts to denote the marginal distributions of $G$.
For example, $G_1$ is the
distribution of $T_1$ and $G_{2,3}$ is the distribution of
$(T_2,T_3)$. For current status censoring ($k=1$), we denote the
observation time simply by $T$.

 We study the MLE $\hat F_n$ of $F_0$, based on $n$ i.i.d.\
copies $W_1,\dots,W_n$ of $W$, where
$W_i = (\mathbf{T}_i,\mathbf{\Delta}_i,Z_i)$, $\mathbf{T}_i =
(T_{1i},\dots,T_{ki})$ and $\mathbf{\Delta}_i =
(\Delta_{1i},\dots,\Delta_{k+1,i})$.
We allow ties between the observation times of $\mathbf{T}_i$ and $\mathbf{T}_j$ for $i\neq j$.

The observed sets $A$ in this model are
given in equations \eqref{eq: observed set Delta+ = 1}
and \eqref{eq: observed set Delta+ = 0}. Recall that $A$ is a line segment
if $\Delta_+=1$ and that $A$ is a half
plane if $\Delta_+=0$. Assuming that $F$ has a density $f$ with
respect to some dominating measure $\mu_X \times \mu_Y$, the
likelihood (up to multiplicative terms only depending on $G$) is given by $L_n(F) = \prod_{i=1}^n
q(W_i)$, where\vspace{-.2cm}
\begin{align*}
   q(w) = q(t,\delta,z) = \prod_{j=1}^k \left\{ \int_{(t_{j-1},t_j]}
   f(s,z)\mu_X(ds) \right\}^{\delta_j}
   \left(1-F_X(t_k)\right)^{1-\delta_{+}},
\end{align*}
and $F_X(x) = F(x,\infty)$ is the marginal distribution of $X$
under $F$. Since the first term of $q$ is a density-type term,
$L_n(F)$ can be made arbitrarily large by letting $f(s,z)$ peak at
$z=Z_i$ for some observation with $\Delta_{+i}=1$. We therefore
define the MLE $\hat F_n(x,y)$ to be the maximizer of $L_n(F)$ over
the class $\mathcal F$ of all bivariate distribution functions that
have a marginal density $f_Y$ with respect to counting measure on
the observed marks. We can then write $L_n(F) = \prod_{i=1}^n
P_F(A_i)$.

Analogously to \citet*{Maathuis05}, we call the projection of $A$ on the $x$-axis
the \emph{$x$-interval} of $A$.
We denote the left endpoint and right endpoint of the $x$-interval by $L$ and
$R$:\vspace{-.15cm}
\begin{align*}
   L & = \sum_{j=1}^{k+1} \Delta_j T_{j-1}, \qquad  R= \sum_{j=1}^{k+1} \Delta_j
   T_j.
\end{align*}
Furthermore, we define a new variable $U$ that will play an important role
in our analysis:\vspace{-.15cm}
\begin{align}\label{eq: def U}
   U = \Delta_+ R + \Delta_{k+1}L.
\end{align}
Let $U_{(1)},\dots,U_{(n)}$
be the order statistics of $U_1,\dots,U_n$ and let
$\mathbf{\Delta}_{(i)}=(\Delta_{1(i)},\dots,\Delta_{k+1,(i)})$,
$Z_{(i)}$, $A_{(i)}$, $L_{(i)}$ and $R_{(i)}$ be the corresponding
values of $\mathbf{\Delta}$, $Z$, $A$, $L$ and $R$. We break ties in
$U$ arbitrarily after ensuring that observations with $\Delta_+=1$
are ordered before those with $\Delta_+=0$. Recall that the maximal
intersections are the local maximum regions of the height map $h: \R^2 \mapsto \N$ of the
observed sets. Since $Y$ is continuous, the observed sets
$A_{(i)}$ with $i\in\mathcal I = \{i\in
\{1,\dots,n\} \st \Delta_{+(i)}=1\}$ are completely distinct with
probability one. Hence, each such $A_{(i)}$ contains exactly one
maximal intersection $M_{(i)}$ of the form:\vspace{-.1cm}
\begin{align}
  \begin{array}{ll}
   M_{(i)} & = (D_{(i)},R_{(i)}] \times \{Z_{(i)}\}, \quad \text{where}\\
   D_{(i)} & = \max\{ \{ L_{(j)} \st j\notin \mathcal I, j<i\} \cup\{ L_{(i)}\}\}.
  \end{array}\label{eq: maximal intersection}
\end{align}
To understand this expression, let $S_{(i)}$ be the collection of
observed sets $A_{(j)}$ with $\Delta_{+(j)}=0$ and $L_{(i)} <
L_{(j)} < R_{(i)}$. If $S_{(i)}=\emptyset$, then the height map is constant on $A_{(i)}$,
and the complete set $A_{(i)}$ is a local maximum region. Hence, in this case $M_{(i)}=A_{(i)}$ and $D_{(i)}=L_{(i)}$.
On the other hand, if $S_{(i)} \neq \emptyset$, then the height map
is increasing on $A_{(i)}$ in the $x$-direction. Hence, in this case $M_{(i)} \subsetneq A_{(i)}$ and
the left endpoint of $M_{(i)}$ is $\max\{L_{(j)}\st A_{(j)}\in S_{(i)}\}$, which
equals $\max \{L_{(j)} \st j\notin \mathcal I, j<i\}$. Note that the right endpoints of
$M_{(i)}$ and $A_{(i)}$ are always identical. Moreover, note that the equations in \eqref{eq:
maximal intersection} imply that the maximal intersections
can be computed in $O(n\log n)$ time, since the most computationally intensive step
consists of sorting the data. This is faster than the
height map algorithm of \citet*{Maathuis05}, due to the special
structure in the data.

Analogously to the situation for univariate right censored data, there is
an extra maximal intersection
$M_{(n+1)} = A_{(n)} = (U_{(n)},\infty)\times \R$ if and only if
$\Delta_{+(n)}=0$. Let $\overline{\mathcal I}$ be the collection of
indices of all maximal intersections. Thus, $\overline{\mathcal I} =
\mathcal I$ if $\Delta_{+(n)}=1$ and $\overline{\mathcal I} =
\mathcal I \cup \{n+1\}$ if $\Delta_{+(n)}=0$. Let $p_i$ be the
probability mass of maximal intersection $M_{(i)}$, $i \in
\overline{\mathcal I}$. Then the likelihood can be written as\vspace{-.2cm}
\begin{align}
  \prod_{i=1}^n P(A_i) = \prod_{i=1}^n  \left(\sum_{j\in
   \overline{\mathcal I}} p_j\ind{M_{(j)}\subseteq A_{(i)}}\right) =
   \prod_{i=1}^n p_i^{\Delta_{+(i)}} \left(  \sum_{j\ge i+1, j\in \overline{\mathcal I}}
   p_j\right)^{1-\Delta_{+(i)}},\label{eq: likelihood continuous mark}
\end{align}
where the second equality follows from the fact that the data are ordered
with respect to the variable
$U$ which was defined in \eqref{eq: def U}. The MLE $\hat p$ maximizes this
expression under the constraints \eqref{eq: constraints}. From the
analogy with the likelihood \eqref{eq: likelihood hood univariate right censored}
it follows immediately that\vspace{-.2cm}
\begin{align}\label{eq: hat pi continuous mark}
   \hat p_i & = \prod_{j=1}^{i-1} \left(1-\frac{\Delta_{+(j)}}{n-j+1}\right)
   \frac{\Delta_{+(i)}}{n-i+1},\qquad i \in \mathcal I,
\end{align}
and $\hat p_{n+1} = 1-\sum_{i \in \mathcal I} \hat p_i$ if
$\Delta_{+(n)}=0$. Equivalently, we can write\vspace{-.2cm}
\begin{align}
   \sum_{j\ge i, j\in \overline{\mathcal I}} \hat p_j = \prod_{j\le i-1} \left( 1-
   \frac{\Delta_{+(j)}}{n-j+1}\right), \qquad
   i \in \overline{\mathcal I}.\label{eq: hat tail}
\end{align}
These formulas are different from (but equivalent to) the ones given
in Section 3.1 of HMG. The form given here has several advantages.
First, the tail probabilities \eqref{eq: hat tail} can be computed
in time complexity $O(n\log n)$, since sorting the data is the most computationally
intensive step. Furthermore, the
current form provides additional insights about the behavior of the
MLE. In particular, it shows that the MLE can be viewed as a right
endpoint imputation estimator (see Remark \ref{remark: imputation
estimator}), and it allows for a derivation of the
almost sure limit of the MLE (see Section~\ref{sec:limit}).

The vector $\hat p$ is uniquely determined. This was noted by
HMG and also follows from our derivation here. We obtain $\hat
F_n(x,y)$ by summing all probability mass of $\hat p$ that falls in the region $(0,x]\times
(-\infty,y]$. We define a marginal MLE for the distribution of $X$ by letting $\hat
F_{Xn}(x) = \hat F_n(x,\infty)$. The estimators $\hat F_n$ and $\hat
F_{Xn}$ can suffer considerably from representational
non-uniqueness, since the maximal intersections $\{ M_{(i)} \st
i\in\mathcal I\}$ are line segments, and the potential maximal intersection $M_{(n+1)}$ is a half plane.
We let $\hat F_n^{\ell}$ denote the estimator that assigns all
mass to the upper right corners of the maximal intersections, since
it is a lower bound for the MLE. Similarly,
we let $\hat F_n^{u}$ denote the estimator that assigns all mass to the lower left
corners of the maximal intersections, since it is an
upper bound for the MLE. The formulas for $\hat F_n^\ell$ and $\hat
F_{Xn}^{\ell}$ can be written as follows:\vspace{-.2cm}
\begin{align}
   1-\hat F^\ell_{Xn}(x) & = \prod_{U_{(i)}\le x} \left( 1 -
   \frac{\Delta_{+(i)}}{n-i+1}\right),\label{eq: 1-hatFXn(x)}\\
   \hat F^\ell_n(x,y) & = \sum_{i=1}^n \hat p_i \ind{U_{(i)} \le x,
   Z_{(i)}\le y} \notag \\
    & = \sum_{U_{(i)}\le x}\prod_{U_{(j)}<U_{(i)}}
    \left(1-\frac{\Delta_{+(j)}}{n-j+1}\right)\frac{\Delta_{+(i)}\ind{Z_{(i)}\le
    y}}{n-i+1}\label{eq: hatFn(x,y)},
\end{align}
using \eqref{eq: hat pi continuous mark}, \eqref{eq: hat tail} and
the definition of $U$ in \eqref{eq: def U}.

\begin{remark}\label{remark: imputation estimator}{\rm
  The MLE $\hat F_n^\ell$ can be viewed as a right endpoint imputation estimator.
  To see this, consider creating a new collection of observed sets $A'_{(i)}$:\vspace{-.15cm}
  \begin{align*}
     A'_{(i)} & = \left\{ \begin{array}{lll}
        \{U_{(i)}\} \times \{Z_{(i)}\} & \text{if} & i\in \mathcal I,\\
        A_{(i)} & \text{if} & i\notin \mathcal I.
     \end{array}\right.
  \end{align*}
  That is, for each $i=1,\dots,n$, we replace $A_{(i)}$ by its right endpoint if $\Delta_{+(i)}=1$, while
  we leave it unchanged if $\Delta_{+(i)}=0$.
  The intersection structures of $\{A_{(i)}\}_{i=1}^n$ and $\{A'_{(i)}\}_{i=1}^n$ are
  identical, meaning that $A_{(i)}\cap A_{(j)} =\emptyset$ if and
  only if $A'_{(i)}\cap A'_{(j)} = \emptyset$, for all $i,j \in
  \{1,\dots,n\}$. Furthermore, the maximal intersections
  of $\{A'_{(i)}\}_{i=1}^n$ are $\{M'_{(i)} = A'_{(i)} \st i\in \overline{\mathcal
  I}\}$. Hence, writing the likelihood
  for the imputed data in terms of $p$ yields exactly the
  same likelihood as \eqref{eq: likelihood continuous mark}.
  This implies that the maximizing vector $\hat p'$ is identical to the
  vector $\hat p$ for the original data. Moreover, the upper right corners
  of $\{M_{(i)}\}$, $i\in \overline{\mathcal I}$
  and $\{M'_{(i)}\}$, $i\in \overline{\mathcal I}$ are identical. Since $\hat F_n^\ell$ assigns all
  mass to the upper right corners of the maximal intersections, it follows that $\hat F_n^{\ell}$ is
  completely equivalent to the MLE for the modified data.
  Finally, note that the right endpoint imputation scheme imputes
  an $x$-value that is
  always at least as large as the unobserved $X$. This explains why the MLE $\hat F_{Xn}^\ell$
  tends to have a negative bias.
}\end{remark}

\section{Inconsistency of the MLE}\label{sec:limit}

In this section we derive necessary and sufficient conditions for consistency
of the MLEs $\hat F_{Xn}^{\ell}$ and $\hat F_n^{\ell}$ (Theorem \ref{th: iff
conditions}). These conditions force a relation between the
unknown distribution $F_0$ and the observation time distribution
$G$. Since such a relation will typically not hold, it follows that
$\hat F_n^{\ell}$ is inconsistent in general. Corollary \ref{cor:
inconsistency current status} further strengthens this result when
$X$ is subject to current status censoring, and shows that in that
case $\hat F_{Xn}^{\ell}$ is inconsistent for any continuous
choice of $F_0$ and $G$. Corollary \ref{cor: bias decreases with k}
shows that the asymptotic biases of $\hat F_{Xn}^{\ell}$ and $\hat
F_n^{\ell}$ converge to zero as the number $k$ of observation times per
subject increases, at least for one particular distribution of $T_1,\dots,T_k$.

The results in this section are based on deriving the limits $F_{X\infty}^{\ell}$ and $F_{\infty}^{\ell}$
for the lower bounds $\hat F_{Xn}^{\ell}$ and $\hat F_n^{\ell}$ of the MLE.
The reason for looking at these lower bounds is that
$\hat F_{Xn}^{\ell}$ and $\hat F_n^{\ell}$ can be expressed in simple closed forms (see
\eqref{eq: 1-hatFXn(x)} and \eqref{eq: hatFn(x,y)}). Moreover, in
many cases representational non-uniqueness disappears in the limit,
so that the limits of $\hat F_{Xn}$ and $\hat F_n$ are unique and
equal to $F_{X\infty}^{\ell}$ and $F_{\infty}^{\ell}$. Necessary and sufficient
conditions for uniqueness of the limit are: (i) all maximal
intersections $M_{(i)}$, $i\in \mathcal I$, converge to points, and (ii)
$\sum_{i\in \mathcal I}\hat p_i \to 1$ as $n \to \infty$. These conditions
are satisfied in Examples 1 and 2 in Section \ref{sec: examples}.
If these conditions
fail, then the upper bounds $F_{X\infty}^u$ and $F_{\infty}^u$ can
be obtained from their lower bounds by reassigning mass from the
upper right corners of the maximal intersections to the lower left
corners. This occurs in Examples 3 and 4 in Section~\ref{sec: examples},
and further details can be found in
\citet*[Section 9.4]{Maathuis06thesis}.

In order to derive $F_{X\infty}^{\ell}$ and $F_{\infty}^{\ell}$ we
start by rewriting equations \eqref{eq: 1-hatFXn(x)} and \eqref{eq:
hatFn(x,y)} in terms of stochastic processes. We introduce the
following notation:\vspace{-.2cm}
\begin{align}\label{eq: def Hn Vn VXn}
   \begin{array}{rll}
      \H_n(x) & = \P_n \ind{U\le x}, & \qquad x \ge 0,\\
      \V_n(x,y) & = \P_n \Delta_+
          \ind{U\le x, Z\le y}, & \qquad x\ge 0, y\in \R,\\
      \V_{Xn}(x) & \equiv \V_n(x,\infty) = \P_n \Delta_+\ind{U\le
      x}, & \qquad x\ge 0,
   \end{array}
\end{align}
where $U$ is defined in \eqref{eq: def U} and $\P_n f(X) =
n^{-1}\sum_{i=1}^n f(X_i)$. Furthermore, let \vspace{-.2cm}
\begin{align}\label{eq: def hat Lambda}
   \widehat \Lambda_n(x,y) = \int_{[0,x]}
   \frac{\V_n(ds,y)}{1-\H_n(s-)} \quad \text{and} \quad
   \widehat \Lambda_{Xn}(x) \equiv \widehat \Lambda_n(x,\infty) =
   \int_{[0,x]}\frac{\V_{Xn}(ds)}{1-\H_n(s-)}.
\end{align}
Since\vspace{-.1cm}
\begin{align*}
   \widehat \Lambda_n(dx,y) = \frac{\P_n \Delta_+ \ind{U=x,Z\le y}}{\P_n\ind{U\ge
   x}} \quad \text{and} \quad \widehat \Lambda_{Xn}(dx)= \frac{\P_n \Delta_+ \ind{U=x}}{\P_n\ind{U\ge
   x}},
\end{align*}
we can write equations \eqref{eq: 1-hatFXn(x)} and \eqref{eq:
hatFn(x,y)} in terms of $\widehat \Lambda_{Xn}$ and $\widehat
\Lambda_n$:\vspace{-.1cm}
\begin{align}
   1-\hat F^\ell_{Xn}(x) & = \prod_{s \le x} \{ 1 -
   \widehat \Lambda_{Xn}(ds) \}, \label{eq: 1-hatFn(x) integral form}\\
   \hat F^\ell_n(x,y) & = \int_{s\le x}\prod_{u<s}
    \{1-\widehat \Lambda_{Xn}(du)\}\widehat \Lambda_n(ds,y) \label{eq: hatFn(x,y) integral form}.
\end{align}
Note that \eqref{eq: 1-hatFn(x) integral form} is analogous to
the Kaplan-Meier estimator for right censored data, and that
\eqref{eq: hatFn(x,y) integral form} is analogous to equation (3.3)
of \citet*{HuangLouis98}. However, our functions $\widehat
\Lambda_{Xn}$ and $\widehat \Lambda_n$ are defined differently,
since they are based on the variable $U$. This difference lies at the root of the
inconsistency problems of the MLE.

The limits of the processes $\H_n$, $\V_n$, $\V_{Xn}$, $\hat
\Lambda_n$, $\hat \Lambda_{Xn}$, $\hat F_n^{\ell}$ and $\hat
F_{Xn}^{\ell}$ are given in the Appendix (Lemmas \ref{lemma: GC} -
\ref{lem: convergence F}) and are denoted by $H$, $V$, $V_X$,
$\Lambda_{\infty}$, $\Lambda_{X\infty}$, $F_{\infty}^{\ell}$ and
$F_{X\infty}^{\ell}$, respectively. Corollaries \ref{cor: F in terms
of Lambda} - \ref{cor: F when current status censoring} in the Appendix
provide various alternative ways to express $F_{\infty}^{\ell}$.

We are now ready to give necessary and sufficient conditions for consistency
of $\hat F_{Xn}^{\ell}$ and $\hat F_n^{\ell}$, after introducing the following notation:\vspace{-.1cm}
\begin{align}
     H(x) & = V_X (x) + \int_{[0,x]} \{1-F_{0X} (s)\} d G_k(s),\label{eq: H}\\
     V(dx,y) & = \sum_{j=1}^k F_0(x,y)d G_j(x) - \sum_{j=2}^k \int_{[0,x]} F_0(s,y) d G_{j-1,j}(s,x), \label{eq: V(ds,y)}\\
     V_X(dx) & = \sum_{j=1}^k F_{0X}(x)d G_j(x) - \sum_{j=2}^k \int_{[0,x]} F_{0X}(s) d G_{j-1,j}(s,x), \label{eq: V1(ds)}
\end{align}
see equations \eqref{eq: V} - \eqref{eq: H appendix} in the Appendix.
Moreover, throughout this section we let $\tau$ be such that
$H(\tau)<1$, we define $0/0=0$ and $f(x-)=\lim_{t\uparrow x} f(t)$
for any function $f: \R \mapsto \R$.

\begin{theorem}\label{th: iff conditions}
   The MLE is inconsistent in general. The MLE $\hat F_{Xn}^\ell$ is consistent
   for $F_{0X}$ on $(0,\tau]$ if and only if the following condition
   holds for all $x\in (0,\tau]$:\vspace{-.1cm}
   \begin{align}
      \Lambda_{X\infty}(x) \equiv \int_{[0,x]} \frac{V_X(ds)}{1-H(s-)} & = \int_{[0,x]}
      \frac{F_{0X}(ds)}{1-F_{0X}(s-)} \equiv \Lambda_{0X}(x).  \label{eq: CondForEquality Fx}
   \end{align}
   The MLE $\hat F_n^\ell$ is consistent for $F_0$ on $(0,\tau]\times\R$ if and only if
   the following condition holds for all $x\in(0,\tau]$, $y\in \R$:\vspace{-.1cm}
   \begin{align}
      \Lambda_{\infty}(x,y) \equiv \int_{[0,x]} \frac{V(ds,y)}{1-H(s-)} & = \int_{[0,x]}
      \frac{F_0(ds,y)}{1-F_{0X}(s-)} \equiv \Lambda_0(x,y)
        \label{eq: CondForEquality F}.
   \end{align}
   Finally, let $x_0 \in (0,\tau]$ with $F_{X\infty}(x_0)>0$. Then
   $\hat F_n^\ell(x_0, y) / \hat F^\ell_{Xn}(x_0)$ is consistent for
   $F_{0Y}(y)$ if $X$ and $Y$ are independent.
\end{theorem}
\begin{proof}
   The one-to-one correspondence between a univariate distribution function and its
   cumulative hazard function implies that $\hat F_{Xn}^\ell$ is
   consistent for $F_{0X}$ if and only if $\Lambda_{X\infty}$ (equation \eqref{eq: Lambda X infty} in the Appendix)
   equals
   the cumulative hazard function $\Lambda_{0X}$ of $F_{0X}$. This gives condition \eqref{eq: CondForEquality Fx}.
   Similarly, it follows that $\hat F_n^\ell(x,y)$ is consistent for
   $F_0(x,y)$ if and only if $\Lambda_{\infty}$ (equation \eqref{eq: Lambda infty} in the Appendix) equals the mark
   specific cumulative hazard function $\Lambda_0$ of $F_0$. This gives condition \eqref{eq: CondForEquality F}.
   The final claim of the theorem follows from equation \eqref{eq: hat F infty
   for independence} in the Appendix. \qed
\end{proof}

Note that conditions \eqref{eq: CondForEquality Fx} and
\eqref{eq: CondForEquality F} are difficult to interpret, since
$F_{0X}$ and $F_0$ enter on both sides of the equations when we plug
in expressions \eqref{eq: H} -- \eqref{eq:
V1(ds)} for $H(s-)$, $V(ds,y)$ and $V_X(ds)$. However, it is clear
that the conditions force a relation between the unknown distribution $F_0$
and the observation time distribution $G$. Such
a relation will typically not hold and cannot be assumed since
$F_0$ is unknown. Hence, it follows that the MLE is inconsistent in general.
The following corollary further strengthens this
result when $X$ is subject to current status censoring.
\begin{corollary}\label{cor: inconsistency current status}
   Let $X$ be subject to current status censoring, and let $F_{0X}$
   and $G$ be continuous. Then
   the MLE $\hat F_{Xn}^\ell$ is inconsistent for any choice of $F_{0X}$ and $G$.
\end{corollary}
\begin{proof}
   Let $\gamma = \inf\{x: F_{0X}(x)>0\} < \tau$. Since $X$ is subject to current status
   censoring and since the distributions $G$ and $F_{0X}$ are continuous,
   condition \eqref{eq: CondForEquality Fx} can be rewritten as
   \begin{align*}
        \int_{(\gamma,x]} \frac{dG(s)}{1-G(s)} = \int_{(\gamma,x]}
        \frac{dF_{0X}(s)}{F_{0X}(s)\{1-F_{0X}(s)\}}, \qquad x\in (\gamma,\tau].
   \end{align*}
   This integral equation is solved by
   \begin{align*}
      -\log\{1-G(x)\} + C = \log\left\{ \frac{F_{0X}(x)}{1-F_{0X}(x)}
      \right\}, \qquad x\in (\gamma,\tau].
   \end{align*}
   This yields $F_{0X}(x) = [1+\exp(-C)\{1-G(x)\}]^{-1}$ for $x\in (\gamma,\tau]$.
   Since there is no finite $C$ such that $F_{0X}(\gamma ) = 0$ holds, it follows that
   condition \eqref{eq: CondForEquality Fx} fails
   for all continuous distributions $G$ and $F_{0X}$. \qed
\end{proof}
Finally, we show that the asymptotic bias of the MLE converges to zero as
the number $k$ of observation times per subject increases, for at
least one particular distribution of $\mathbf{T} =
(T_1,\dots,T_{k})$, namely if $T_1,\dots,T_k$ are distributed as
the order statistics of a uniform sample on $[0, \theta]$. The proof
of this result is given in the Appendix.
\begin{corollary}\label{cor: bias decreases with k}
  Let $X$ be subject to interval censoring case $k$, and let the elements
  $T_1,\dots,T_k$ of $\mathbf{T}$ be the order statistics of $k$
  independent uniform random variables on $[0,\theta]$.
  Let $V^k(x,y)$, $V^k_X(x)$, $H^k(x)$, $\Lambda_{\infty}^k(x,y)$ and
  $\Lambda_{X\infty}^k(x)$ denote the limits defined in Lemmas \ref{lemma: GC} and
  \ref{lem: convergence lambda}, using the superscript $k$ to denote the
  dependence on $k$. Then
  \begin{align*}
    \Lambda_{X\infty}^k (x)
    = \int_{[0,x]} \frac{ V_X^k (ds) }{1-H^k(s-) } & \rightarrow \int_{[0,x]}  \frac{ F_{0X}(ds)}{1-F_{0X}(s-)} =
    \Lambda_{0X}
    (x), \quad k\to\infty,\\
    \Lambda_{\infty}^k (x,y) =
    \int_{[0,x]} \frac{ V^k (ds,y) }{1-H^k(s-) } & \rightarrow
        \int_{[0,x]} \frac{F_0(ds,y)}{1-F_{0X}(s-)} = \Lambda_0(x,y),
        \quad k\to\infty,
  \end{align*}
  for all continuity points of $\Lambda_{0X}$ and $\Lambda_{0}$ with $x<\theta$ and $y\in \R$.
\end{corollary}

\section{Repaired MLE via discretization of marks}\label{sec: repaired mle}

We now define a simple repaired estimator $\tilde F_n(x,y)$ which is
consistent for $F_0(x,y)$ for $y$ on a grid. The idea behind the
estimator is that one can define discrete competing risks based on a
continuous random variable. Doing so transforms interval censored
continuous mark data into interval censored data with competing risks.

To describe the method, we let $K>0$ and define a grid $-\infty
\equiv y_0 < y_1 < \dots < y_K < y_{K+1} \equiv \infty$. Next, we
introduce a new random variable $C\in\{1,\dots,K+1\}$:
\begin{align*}
   C = \sum_{j=1}^{K+1} j \ind{y_{j-1}<Y \le y_j}.
\end{align*}
We can determine the value of $C$ for all observations with an
observed mark. Hence, we can transform the observations
$(\mathbf{T},\mathbf{\Delta},Z)$ into
$(\mathbf{T},\mathbf{\Delta},Z^*)$, where $Z^* = \Delta_+ C$. This
gives interval censored data with $K+1$ competing risks.

Since the observed sets for interval censored data with competing risks
form a partition of the space $\Rplus \times \{1,\dots,K+1\}$,
Hellinger consistency of the MLE follows from Theorems 9 and 10 of
\citet*{VanderVaartWellner00:GC-preservationTheorems}. Under some
additional regularity conditions, we can derive local and uniform
consistency from the Hellinger consistency, see \citet*[Section 4.2]{Maathuis06thesis}.
This means that we can consistently estimate the sub-distribution
functions $F_{0j}(x) = P(X\le x, C=j) = P(X\le x, y_{j-1}<Y\le
y_j)$, $x\in \Rplus$. Hence, we can consistently estimate $F_0(x,y_j) =
\sum_{\ell=1}^j F_{0\ell}(x)$ for $x\in \Rplus$ and $y_j$ on the
grid.

Note that the introduction of the variable $C$ causes more
overlap between observed sets, since previously non-overlapping horizontal line
segments may overlap if they are assigned the same value of $C$.
As a result, the repaired MLE has
smaller maximal intersections in the $x$-direction. Hence, the repaired
MLE is affected less by representational non-uniqueness on the $x$-axis. This
is visible in Examples 3 and 4 in Section~\ref{sec: examples}.

The repaired MLE can be computed with one of the algorithms
described in \citet*[Section 2.4]{GroeneboomMaathuisWellner06a}. It
may be tempting to choose $K$ large, such that $F_0(x,y)$ can be
estimated for $y$ on a fine grid. However, this may result in a poor
estimator. To obtain a good estimator one should choose the grid
such that there are ample observations for each value of $C$. In
practice, one can start with a coarse grid, and then refine the grid
as long as the estimator stays close to the one computed on the
coarse grid.

In principle it is possible to estimate the entire joint
distribution function $F_0(x,y)$ for $(x,y)$ in the interior of the
support of the distribution of the observation times under smoothness
assumptions on $F_0$. This would proceed by letting both $K$ and the
$y_j$'s defining the partition all depend on $n$ in such a way that
$K=K_n \to \infty$,
\begin{align*}
  \max_{1 \le j \le K_n-1}
  (y_{j+1,n}- y_{j,n} ) & \rightarrow 0, \qquad \text{and} \qquad
  n \min_{1 \le j \le K_{n}-1} (y_{j+1,n} - y_{j,n} ) \rightarrow
  \infty,
\end{align*}
as $n\to \infty$. It would even be possible to choose $K_n$ and
$\{y_{j,n} \}$ depending on the data via model-selection methods
(see, e.g., \citet*{BirgeMassart97} and
\citet*{BarronBirgeMassart99}), but these further developments are
beyond the scope of the present paper and will be investigated in
detail elsewhere.

\cite*{Maathuis06thesis},
\citet*{GroeneboomMaathuisWellner06a} and \citet*{GroeneboomMaathuisWellner06b}
showed that the MLE for current status data with competing risks
converges at rate $n^{1/3}$ to a new self-induced limiting distribution.
This result implies that one can use subsampling
to construct pointwise confidence intervals for the
sub-distribution functions (\citet*{PolitisRomanoWolf99}). This method is also valid
for the repaired MLE for current status data with continuous marks, and can be used for the
construction of pointwise confidence intervals for $F_0(x,y)$
for $y$ on the grid. The limiting distribution of the MLE for more general forms
of interval censoring with competing risks has not yet been established, and in such
cases the use of subsampling is therefore
not yet justified.

\citet*{JewellVanderLaanHenneman03} and \citet*[Chapter
7]{Maathuis06thesis} studied estimation of a family of smooth
functionals of the sub-distribution functions for current status
data with competing risks. \citet{JewellVanderLaanHenneman03}
suggested that their ``naive estimator" yields asymptotically
efficient estimators for these smooth functionals, and
\citet*{Maathuis06thesis} showed that the same is true for the MLE.
These results extend to the repaired MLE for current status data
with continuous marks. Asymptotic properties of estimators of smooth
functionals for more general forms of interval censoring with
competing risks are currently still unknown.

\section{Examples}\label{sec: examples}

In this section we support the theoretical results of Sections
\ref{sec:limit} and \ref{sec: repaired mle} by simulations. In
particular, we show support for our claims that $\hat F_n^{\ell}
\convas F_{\infty}^{\ell}$, $\hat F_n^{u} \convas F_{\infty}^{u}$
and $\tilde F_n \convas F_0$. Moreover, we show that the difference
between the true underlying distribution $F_0$ and the limits of the
MLE $F_{\infty}^{\ell}$ and $F_{\infty}^{u}$ can be considerable. We
give four examples that cover a wide range of scenarios. They
include cases where $X$ and $Y$ are independent (Ex. 1) or dependent
(Ex. 2 -- 4), where $X$ is subject to interval censoring case 1 (Ex.
1, 2) or case 2 (Ex. 3, 4), and where the distribution of ${\bf T}$
is continuous (Ex. 1 -- 3) or discrete (Ex. 4).

\begin{example}\label{ex: example 1}{\rm
  Let $X$ and $Y$ be independent, with $X\sim \text{Unif}(0,1)$ and $Y\sim\text{Exp}(1)$.
  Let $X$ be subject to current status censoring with observation
  time $T \sim \text{Unif}(0,0.5)$ independent of $(X,Y)$.
}\end{example}

\begin{example}\label{ex: example 2}{\rm
   Let $X \sim \text{Unif}(0,1)$, and let $Y|X$ be exponentially distributed
   with mean $2/(2X+1)$.
   Let $X$ be subject to current status censoring with observation time
   $T \sim \text{Unif}(0,1)$ independent of $(X,Y)$.
}\end{example}

\begin{example}\label{ex: example 3}{\rm
   Let $X\sim \text{Unif}(0,2)$, and let $Y\equiv X$. Let $X$ be subject to interval censoring
   case 2 with observation times $(T_1,T_2)$,
   independent of $(X,Y)$ and uniformly distributed over
    $\{(t_1,t_2)\st 0\le t_1\le 1, 1\le t_2 \le 2\}$.
}\end{example}

\begin{example}\label{ex: example 4}{\rm
   Let $(X,Y)$ be uniformly distributed over $\{(x,y) \st 0\le x \le y \le 1\}$.
   Let $X$ be subject to interval censoring case 2 with observation times
   $(T_1,T_2)$ independent of $(X,Y)$.
   Let the distribution of $(T_1,T_2)$ be discrete: $G\{(0.25, 0.5)\} =0.3$,
   $G\{(0.25, 0.75)\} = 0.3$ and $G\{(0.5,
   0.75)\} = 0.4$.
}\end{example}

\begin{figure}[ht]
\begin{center}
    \includegraphics[scale=.65]{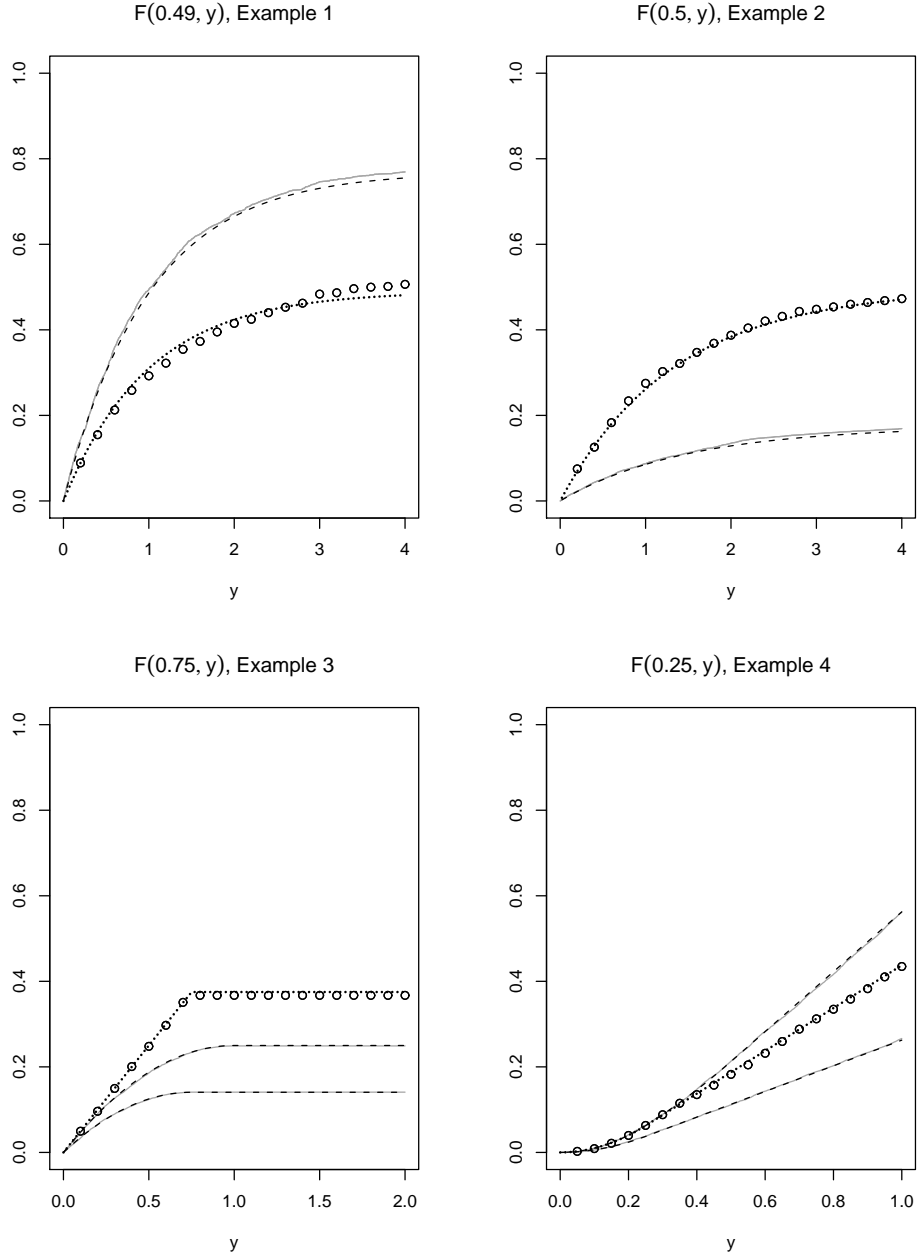}
    \caption{ Estimation of $F_0(x_0,y)$ in Examples 1 -- 4, for fixed $x_0$ and $y\in \R$.
     Dotted: the true underlying
    distribution $F_0(x_0,y)$. Solid grey:
    the MLEs $\hat F^\ell_{n}(x_0,y)$ and $\hat F^u_{n}(x_0,y)$. Dashed: the limits $F^\ell_{\infty}(x_0,y)$ and
    $F^u_{\infty}(x_0,y)$ of the MLE.
    Circles: the repaired MLE $\tilde F^\ell_{n}(x_0,y) =  \tilde F^u_{n}(x_0,y)$,
    using an equidistant grid with
    $K=20$.
    In all cases $n=10,\!000$.}
    \label{fig: Fx0y}
\end{center}
\end{figure}

For each example we derived the limits $F_{\infty}^{\ell}$
and $F_{\infty}^u$ of the MLE, using Lemma \ref{lem: convergence F}. Details of
these derivations are given in \citet*[Section
9.4]{Maathuis06thesis}. We also computed the MLEs $\hat F_n^\ell$
and $\hat F_n^u$ and the repaired MLEs $\tilde F_n^\ell$ and
$\tilde F_n^u$ for a simulated data set of size $n=10,\!000$. For
the repaired MLE we used an equidistant grid with $K=20$ points
as shown in Figure \ref{fig: Fx0y}.

\begin{figure}[ht]
\begin{center}
    \includegraphics[scale=.65]{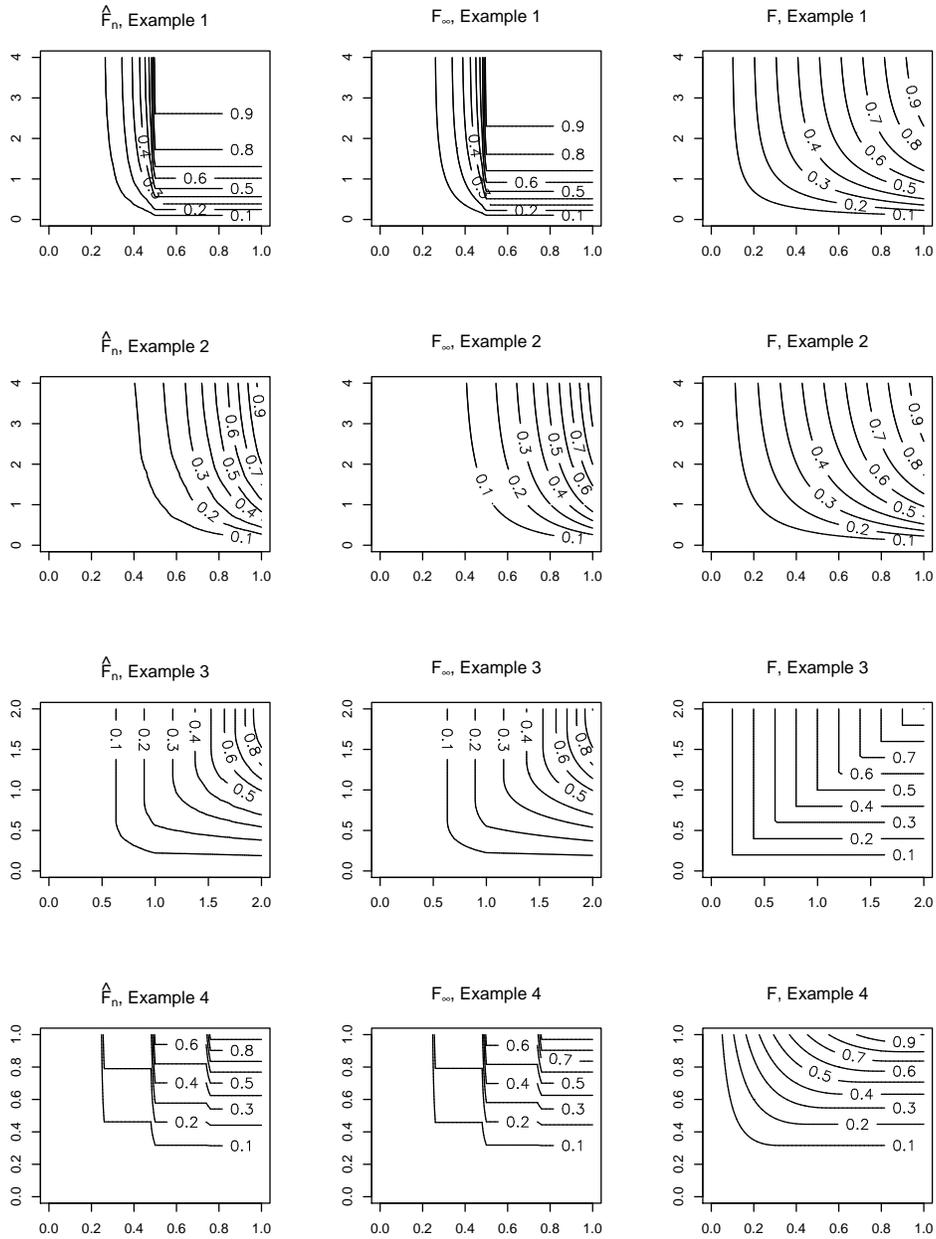}
    \caption{Contour lines of the bivariate functions $\hat F^\ell_n$ (left column), $F^\ell_{\infty}$ (middle column)
    and $F_0$ (right column) for Examples 1 -- 4. All
    functions were computed on an equidistant grid with grid size 0.02, and sample size $n=10,\!000$. }
    \label{fig: F2D}
\end{center}
\end{figure}

\begin{figure}[ht]
\begin{center}
    \includegraphics[scale=.65]{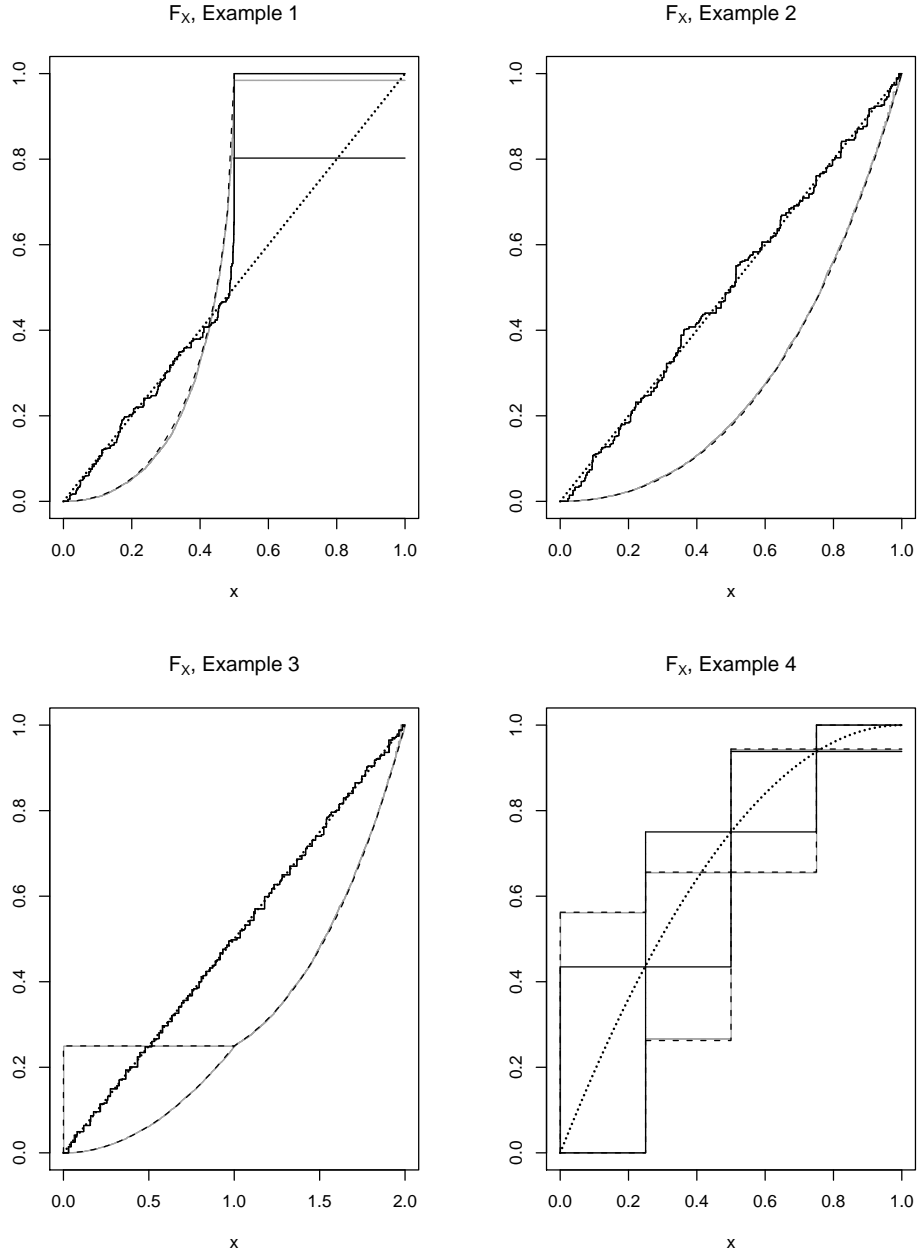}
    \caption{Estimation of $F_{0X}$ in Examples 1 -- 4. Dotted: the
   true underlying
    distribution $F_{0X}$. Solid grey:
    the MLEs $\hat F^\ell_{Xn}$ and $\hat F^u_{Xn}$. Dashed: the limits $F^\ell_{X\infty}$ and
    $F^u_{X\infty}$ of the MLE.
    Solid black: the repaired MLEs $\tilde F^\ell_{Xn}$ and $\tilde F^u_{Xn}$, using the equidistant grid with
    $K=20$ shown in Figure \ref{fig: Fx0y}.
    In all cases $n=10,\!000$.}
    \label{fig: Fx}
\end{center}
\end{figure}

The results are given in Figures
\ref{fig: F2D} - \ref{fig: Fx0y}. These figures show that the MLEs $\hat F_n^{\ell}$
and $\hat F_n^{u}$ are indeed very close to our derived limits
$F_{\infty}^{\ell}$ and $F_{\infty}^{u}$. On the other hand, the
repaired MLEs $\tilde F_n^{\ell}$ and $\tilde F_n^u$ are very close
to the true underlying distribution $F_0$. Moreover, the results show that there
can be a very significant difference between the limit of the MLE and
the true underlying distribution $F_0$.

We now discuss the simulation results in more detail. Figure \ref{fig: F2D}
considers estimation of the joint distribution $F_0$. It shows the contour lines of
the MLE $\hat F_n^\ell$, its limit
$F_{\infty}^\ell$, and the true underlying distribution $F_0$. Note
that $\hat F_n^\ell$ and $F_{\infty}^\ell$ are almost
indistinguishable, while there is a clear difference between
$F_{\infty}^\ell$ and $F_0$. The results for the upper limits $\hat
F_n^u$ and $F_{\infty}^u$ are similar and not shown. Results for the repaired MLE are
not shown since this estimator only takes values for $y$ on a grid.

Figure \ref{fig: Fx} considers estimation of the marginal distribution $F_{0X}$. We see that
the MLEs $\hat F_{Xn}^{\ell}$ and $\hat F_{Xn}^{u}$ are close to the derived limits $F_{X\infty}^{\ell}$ and
$F_{X\infty}^{u}$. Moreover, note that $\hat F_{Xn}^{\ell}$ tends to be below $F_{0X}$.
This can be understood
via Remark \ref{remark: imputation estimator} on page
\pageref{remark: imputation estimator}, which explains that $\hat
F_n^{\ell}$ can be viewed as a right endpoint estimator, and hence tends to have a negative bias.
Note that the repaired MLE $\tilde F_n$ closely follows
$F_{0X}$.

Figure \ref{fig: Fx0y} considers estimation of $F_0(x_0,y)$ for fixed
$x_0$. The function $F_0(x_0,y)$ is often estimated as an alternative for
$F_{0Y}$, since $F_{0Y}$ is heavily affected by representational non-uniqueness if the
support of $T_1,\dots,T_k$ is strictly contained in the support of $X$,
a situation that often occurs in practice. The values of $x_0$
were chosen to show a range of scenarios for the behavior
of the MLE, and we see that $\hat F_n(x_0,y)$ can be much too large, much too small and non-unique.
The repaired MLE $\tilde F_n$ is again close to the underlying distribution.

Note that our examples are not linked to any
specific application. For readers who are interested in a comparison
between the MLE and the repaired MLE in a practical situation, we
refer to HMG. They provide such a comparison for the HIV/AIDS
vaccine trial data VAX004 (\citet{FlynnForthalHarroJudsonMayerPara05}), as
well as for simulated data that mimic the vaccine data. They show a difference
between the MLE and the repaired MLE in this setting, but the size of the
difference is quite small.
This can be explained by Corollary \ref{cor: bias decreases with k}, since the time
between successive follow-up visits is relatively short (about 6 months) and the infection
rate is low. Much
larger differences can be expected in, for example, cross-sectional HIV studies, where there is
only one observation time per person.

{\bf [Figure \ref{fig: F2D} -- \ref{fig: Fx0y} about here.]}

\section{Discussion}\label{sec: discussion}

We studied the MLE of the bivariate distribution
of an interval censored survival time and a continuous
mark variable. We derived the almost sure limit of the MLE, and
showed that the MLE is inconsistent in general. We proposed a
simple method to repair the inconsistency, and illustrated the
behavior of the inconsistent and repaired MLE in four examples.

We were prompted to investigate consistency of the MLE in the interval censored
continuous mark model, since the observed sets in this model can
take the form of line segments. Such line segments are an indicator of consistency problems for the MLE,
since the MLE for bivariate
censored data has been found to be inconsistent before when such line segments were present (\citet*{VanderLaan96}
and \citet*[Section 6.2]{Maathuis03}). In this sense our results do not come as a surprise, and they
confirm the idea that the presence of line segments is indicative of consistency problems of the MLE.

There are, however, interesting differences in the underlying
reasons for inconsistency in the above mentioned models. The
inconsistency of the MLE in the model considered by
\citet*{Maathuis03} could be explained by representational
non-uniqueness of the MLE. This is not the case for the interval
censored continuous mark model, where the MLE is typically
inconsistent even if its limit is fully unique. Rather, the
inconsistency in the interval censored continuous mark model can be
explained by the fact that the cumulative hazard functions that
define the MLE in \eqref{eq: 1-hatFXn(x)} and \eqref{eq: hatFn(x,y)}
do not converge to the true underlying cumulative hazard functions.

Finally, we provide a more detailed discussion of the connections
between the current paper and the paper by HMG, since these papers
have been heavily influenced by each other. HMG started studying the
interval censored continuous mark model, in order to analyze data
from the first Phase III HIV/AIDS vaccine trial VAX004
(\citet{FlynnForthalHarroJudsonMayerPara05}). We suspected
inconsistency of the MLE in this model, and investigated this issue
more closely. This study has resulted in the current paper. In turn,
our paper has influenced the work of HMG and their analysis of the
VAX004 data.

There are also some differences between the models in the two
papers. HMG considered a slightly more complicated interval censored
continuous mark model, assuming that $X$ is mixed case interval
censored (as discussed in Section 1) instead of case $k$ interval
censored. They showed that our results in Sections \ref{sec:limit}
and \ref{sec: repaired mle} can be generalized to that situation.
Thus, the MLE is typically inconsistent in this model as well, and
this inconsistency can be repaired by discretizing the marks. HMG
also considered a complication regarding the mark variable $Y$. In
addition to assuming that $Y$ is missing for all individuals who did
not experience the failure event, they allowed $Y$ to be missing
with some probability $p \in (0,1)$ for individuals who did
experience the failure event. In this case there is no closed form
available for the MLE. It is therefore more difficult to study
consistency issues, and consistency of the MLE in this model is
currently still an open problem. However, due to the presence of
line segments we expect inconsistency, and this conjecture is
supported by simulation results of HMG. HMG therefore included our
repaired MLE in the analysis of the VAX004 data.

\section{Acknowledgements}
This research was supported by NSF grant DMS-0203320. We would like
to thank Piet Groeneboom and Michael Hudgens for helpful discussions
and comments.  We also owe thanks to an anonymous referee and an
Associate Editor for useful suggestions concerning the presentation
of our results.

\vspace{1cm} \noindent M.H. Maathuis, University of Washington,
Department of Statistics, Campus Box 354322, Seattle, WA
98195-4322, USA.\\
Email: marloes@stat.washington.edu

\section{Appendix}\label{sec: appendix}

This section contains several technical lemmas and proofs that are needed for the results in Section
\ref{sec:limit}. Lemma \ref{lemma: GC} gives the almost sure limits $H$, $V$ and $V_X$ of the processes
$\H_n$, $\V_n$, $\V_{Xn}$ that were defined in \eqref{eq: def Hn Vn VXn}. Lemma \ref{lem: convergence lambda}
provides the almost sure limits $\Lambda_{\infty}$ and $\Lambda_{X\infty}$ of the processes $\hat \Lambda_{n}$
and $\hat \Lambda_{Xn}$ that were defined in \eqref{eq: def hat Lambda}.
Lemma \ref{lem: convergence F} gives the almost sure
limits $F_{\infty}^{\ell}$ and $F_{X\infty}^{\ell}$ of the MLEs $\hat F_n^{\ell}$ and
$\hat F_{Xn}^{\ell}$ that were given in \eqref{eq: 1-hatFXn(x)} and \eqref{eq: hatFn(x,y)}.
Corollary \ref{cor: F in terms of Lambda} provides an alternative way to express $F_{\infty}^{\ell}$.
Corollaries \ref{cor: yF when independent} and \ref{cor: F when current status censoring}
specialize this result to two special cases, namely the case that $X$ and $Y$ are independent, and the case
that $X$ is subject to current status censoring. Finally, we provide a proof
of Corollary \ref{cor: bias decreases with k}.

\begin{lemma}\label{lemma: GC}
   For $I\subseteq \R^d$ with $d\ge 1$, and let $\mathcal D(I)$ be the space of
   cadlag functions on $I$.
   Furthermore, let $\|\cdot \|_{\infty}$ be the supremum norm on
   $(\mathcal D(\Rplus),\mathcal D(\Rplus),\mathcal D (\Rplus \times \R))$.
   Then
   \begin{align}
       \| (\H_n - H, \V_{Xn} - V_X, \V_n - V )
      \|_{\infty} \convas 0\label{eq:       GlivenkoCantelli},
   \end{align}
   where
   \begin{align}
     V(x,y) & = \sum_{j=1}^k \int_{[0,x]} F_0(t,y) dG_j(t) - \sum_{j=2}^k \int_{0\le s\le t\le x} F_0(s,y)
         dG_{j-1,j}(s,t), \label{eq: V}\\
      V_X(x) & = \sum_{j=1}^k \int_{[0,x]} F_{0X}(t) dG_j(t) - \sum_{j=2}^k \int_{0\le s
      \le t\le x} F_{0X}(s) d
         G_{j-1,j}(s,t), \label{eq: VsubOne}\\
     H(x) & = V_X (x) + \int_{[0,x]} \{1-F_{0X} (s)\} d G_k
     (s),\label{eq: H appendix}
   \end{align}
   and $G_{j-1,j}$ and $G_k$ are defined in the
   beginning of Section~\ref{subsec:closedformcontinuousmark}.
\end{lemma}
\begin{proof}
   Equation \eqref{eq: GlivenkoCantelli} follows immediately from the
   Glivenko-Cantelli theorem, with $H(x) = E(\ind{U\le x})$,
   $V(x,y) = E(\Delta_+\ind{U\le x,Z\le y})$ and $V_X(x) =
   V(x,\infty) = E(\Delta_+ \ind{U\le x})$. We now express $H$,
   $V$ and $V_X$ in terms of $F_0$ and $G$. Note
   that the events $[\Delta_j =1]$, $j=1,\ldots , k+1$, are
   disjoint. Furthermore, note that $U = T_j $ and $Z=Y$ on $[\Delta_j =1]$, $j = 1, \ldots , k$, and
   $U = T_k$ on $[ \Delta_{k+1} =1]$. Hence,
   \begin{align*}
      V(x,y) & = E (\Delta_+  \ind{ U \le x, Z \le y })
                    =  \sum_{j=1}^k P( \Delta_j = 1, Y \le y, T_j \le x) \\
      & = \sum_{j=1}^k P( X \in (T_{j-1} ,T_j], Y \le y, T_j \le x)\\
      & = \sum_{j=1}^k \int_{0 \le s \le t \le x}
      \{ F_0 (t,y) - F_0(s,y) \} d G_{j-1,j} (s,t).
   \end{align*}
   Using $T_0 = 0$, $X>0$ and $G(\{0 < T_1 < \dots < T_k\})=1$,
   this can be written as
   \begin{align*}
      \sum_{j=1}^k \int_{[0,x]} F_0(t,y) dG_j(t) - \sum_{j=2}^k \int_{0\le s\le t \le x}
      F_0(s,y)
         dG_{j-1,j}(s,t).
   \end{align*}
   Taking $y=\infty$ yields the expression for $V_X(x)$. The expression for $H$
   follows similarly, using
   \begin{align*}
      H(x) & = E \ind{ U \le x }
                  =  \sum_{j=1}^k P( \Delta_j = 1, T_j \le x) + P( \Delta_{k+1} = 1, T_k \le
                  x).
   \end{align*}
   \qed
\end{proof}

\begin{lemma}\label{lem: convergence lambda}
   Let $\|\cdot\|_{\infty}$ be the supremum norm on $(\mathcal
   D[0,\tau],  \mathcal D([0,\tau]\times \R))$. Then
   \begin{align*}
      \| (\widehat \Lambda_{Xn} - \Lambda_{X\infty}, \widehat \Lambda_n - \Lambda_{\infty})
      \|_{\infty} \convas 0,
   \end{align*}
   where
   \begin{align}
       \Lambda_{\infty}(x,y) &  = \int_{[0,x]} \frac{V(ds,y)}{1-H(s-)}, \qquad \quad \qquad \qquad \,\, x\in[0,\tau], y\in \R,
        \label{eq: Lambda infty}\\
       \Lambda_{X\infty}(x) & = \Lambda_{\infty}(x,\infty) = \int_{[0,x]} \frac{V_X(ds)}{1-H(s-)},  \qquad x\in[0,\tau].
       \label{eq: Lambda X infty}
   \end{align}
\end{lemma}
\begin{proof}
   This proof is similar to the discussion on page 1536 of
   \citet*{GillJohansen90}. For all $x\ge 0$, let $\H_n^-(x) \equiv
   \H_n(x-)$.
   Consider the mappings
   \begin{align*}
      \left(\H_n^-, \V_{Xn}, \V_n\right) \to \left(\{1-\H_n^-\}^{-1}, \V_{Xn}, \V_{n}\right) \to
      \left(\widehat \Lambda_{Xn}, \widehat \Lambda_n\right)
   \end{align*}
   on the spaces
   \begin{align*}
      \left(\mathcal D_{-}[0,\tau],\mathcal D[0,\tau], \mathcal
      D([0,\tau]\times \R)\right) & \to \left(\mathcal D_{-}[0,\tau],
      \mathcal D[0,\tau], \mathcal D([0,\tau]\times \R)\right)\\
      & \to \left(\mathcal D[0,\tau], \mathcal
      D([0,\tau]\times \R)\right),
   \end{align*}
   where $\mathcal D_{-}(0,\tau]$ is the space of `caglad' (left-continuous with right limits)
   functions on $(0,\tau]$.
   The first mapping is continuous with respect to the
   supremum norm when we restrict the domain of its first argument to elements of
   $\mathcal D_{-}[0,\tau]$ that are bounded by say $\{1+H(\tau)\}/2  < 1$.
   Strong consistency of $\H_n^-$
   ensures that it satisfies this bound with probability one for $n$ large enough.
   The second mapping is continuous with respect to the supremum
   norm by the Helly-Bray lemma. Combining the continuity of these mappings
   with Lemma \ref{lemma: GC} yields the result of the
   theorem.\qed
\end{proof}

\begin{lemma}\label{lem: convergence F}
   Let $\|\cdot\|_{\infty}$ be the supremum norm on $\mathcal
   (D[0,\tau], \mathcal D([0,\tau]\times \R))$. Then
   \begin{align*}
      \| ( \hat F_{Xn}^\ell - F_{X\infty}^\ell , \hat F_n^\ell - F_{\infty}^\ell )\|_{\infty}
      \convas
      0,
   \end{align*}
   where
   \begin{align}
      F_{X\infty}^\ell(x) &  =  1 - \prod_{s\le x} \left\{ 1 - \Lambda_{X\infty}(ds)\right\}, \label{eq: FXinfty}\\
      F_{\infty}^\ell(x,y) & = \int_{u\le x} \prod_{s < u} \left\{ 1- \Lambda_{X\infty}(ds)\right\}\Lambda_{\infty}(du,y).
      \label{eq: Finfty}
   \end{align}
\end{lemma}
\begin{proof}
   To derive the almost sure limit of $\hat F_{Xn}^{\ell}$, consider the mapping
   \begin{align}
       \widehat \Lambda_{Xn} \to \prod_{s\le x} \{ 1-\widehat
       \Lambda_{Xn}(ds)\} = 1-\hat F^\ell_{Xn}(x) \label{eq: mapping Lambda}
   \end{align}
   on the space $\mathcal D[0,\tau]$ to itself.
   This mapping is continuous with respect to the supremum norm
   when its domain is restricted to functions
   of uniformly bounded variation (\citet*{GillJohansen90}, Theorem 7).
   Note that, for $s\in [0,\tau]$, $\widehat \Lambda_{Xn}(s) \le 1/\{1-\H_n(\tau)\} <
   2/\{1-H(\tau)\}$ with probability one for $n$ large enough.
   Together with the monotonicity of $\widehat \Lambda_{Xn}$ this
   implies that with probability one $\widehat \Lambda_{Xn}$ is of uniformly bounded
   variation on $[0,\tau]$, for $n$ large enough. The almost sure limit of $\hat F_{Xn}^\ell$ now follows by
   combining Lemma \ref{lem: convergence lambda} and
   the continuity of \eqref{eq: mapping Lambda}.

   To derive the almost sure limit of $\hat F_n^\ell$ consider the mapping
   \begin{align*}
      (\widehat \Lambda_{Xn}, \widehat \Lambda_{n}) \to \int_{u\le x} \prod_{s < u}
      \{ 1- \widehat \Lambda_{Xn}(ds)\} \widehat \Lambda_n(du,y) =
      \hat F_n^\ell(x,y)
   \end{align*}
   on the
   space $(\mathcal D[0,\tau], \mathcal D([0,\tau]\times \R))$ to $\mathcal D([0,\tau]\times \R)$.
   This mapping is continuous with respect to the supremum norm
   when its domain is restricted to functions of uniformly bounded
   variation (\citet*{HuangLouis98}, Theorem 1). Note that $\widehat \Lambda_n(x,y) \le
   \widehat \Lambda_{Xn}(x)$, so that with probability one the pair $(\widehat \Lambda_n, \widehat \Lambda_{Xn})$
   is uniformly bounded for $n$ large enough. The result then
   follows as in the first part of the proof.
  \qed
\end{proof}

\begin{corollary}\label{cor: F in terms of Lambda}
   For $x\in [0,\tau], y\in \R$, we can write
   \begin{align}
      F^\ell_{\infty}(x,y) =
      \int_{[0,x]}\frac{\Lambda_{\infty}(ds,y)}{\Lambda_{X\infty}(ds)}
      dF^\ell_{X\infty}(s) =
      \int_{[0,x]}
      \frac{V(ds,y)}{V_X(ds)}d F^\ell_{X\infty}(s). \label{eq: hat F infty in terms of V}
   \end{align}
\end{corollary}
\begin{proof}
   Combining equations \eqref{eq: FXinfty}
   and \eqref{eq: Finfty} yields
   \begin{align}\label{eq: Finfty in terms of FXinfty}
      F^\ell_{\infty} (x,y) & =  \int_{[0,x]} \{1-F^\ell_{X\infty}(s-)\} \Lambda_{\infty}(ds,y) \, .
   \end{align}
   Taking $y = \infty$ gives $F^\ell_{X\infty} (x) =  F^\ell_{\infty} (x, \infty)  =  \int_{[0,x]} \{1 -
    F^\ell_{X\infty}(s-)\} \Lambda_{X\infty}(ds)$, so that
    $dF^\ell_{X\infty}(s) =
    \{1-F^\ell_{X\infty}(s-)\}\Lambda_{X\infty}(ds)$. Combining this
    with equation \eqref{eq: Finfty in terms of FXinfty} yields
    the first equality of \eqref{eq: hat F infty in terms of V}.
    The second equality
    follows from the identities
    \begin{align*}
       \Lambda_{\infty}(ds,y) & = V(ds,y) / \{1-H(s-)\},\\
       \Lambda_{X\infty}(ds) & = V_X(ds,y) / \{1-H(s-)\}.
    \end{align*}
    \qed
\end{proof}

\begin{corollary}\label{cor: yF when independent}
   Let $X$ and $Y$ be independent. Then
   \begin{align}
      F^\ell_{\infty}(x,y) = F^\ell_{X\infty}(x)F_{0Y}(y), \qquad x\in [0,\tau], y\in \R.
      \label{eq: hat F infty for independence}
   \end{align}
\end{corollary}
\begin{proof}
   If $X$ and $Y$ are independent, equations \eqref{eq: V(ds,y)} and \eqref{eq: V1(ds)} yield
   $V(ds,y) = F_{0Y}(y)V_X(ds)$.
   Substituting this into equation \eqref{eq: hat F infty in terms of V}
   gives the result. \qed
\end{proof}

\begin{corollary}\label{cor: F when current status censoring}
   Let $X$ be subject to current status censoring ($k=1$). Then
   \begin{align*}
      F^\ell_{\infty}(x,y) = \int_{[0,x]} P(Y\le y|X\le s)d
      F^\ell_{X\infty}(s), \qquad x\in [0,\tau], y\in \R.
   \end{align*}
\end{corollary}
\begin{proof}
   For $k=1$ equations \eqref{eq: V(ds,y)} and \eqref{eq: V1(ds)} reduce to
   $V(ds,y) =  F_0(s,y)dG(s)$ and $V_X(ds) = F_{0X}(s) dG(s)$.
   Hence, $ V(ds,y) / V_X(ds) = F_0(s,y) / F_{0X}(s)  = P(Y\le y|X\le
   s)$. Substituting this into equation \eqref{eq: hat F infty in terms of V}
   completes the proof.
   \qed
\end{proof}

\begin{proof_of}[Corollary \ref{cor: bias decreases with k}]
   Since the observation times are the order statistics of $k$ i.i.d.\
   uniform random variables, the marginal
   densities $g_j$, $j=1,\dots,k$ and the joint densities $g_{j-1,j}$,
   $j=2,\dots,k$ are known (see, e.g., \citet{ShorackWellner86}, page
   97). Summing them over $j$ yields:
   \begin{align*}
     & \sum_{j=1}^k g_j (t)  =
     \frac{k}{\theta}1_{[0,\theta]} (t)  \sum_{j-1=0}^{k-1} {k-1 \choose j-1}
              \left ( \frac{t}{\theta} \right )^{j-1}   \left (1 - \frac{t}{\theta} \right )^{k-1-(j-1)}
       =   \frac{k}{\theta}1_{[0,\theta]} (t),\\
     & \sum_{j=2}^k g_{j-1,j} (s,t)
         =   \frac{k(k-1)}{\theta^2} 1_{[0 \le s \le t \le \theta]}
            \left ( 1 - \frac{t-s}{\theta} \right )^{k-2} .
   \end{align*}
   Let $x<\theta$. Plugging the above expressions for $g_j$ and
   $g_{j-1,j}$ into \eqref{eq: V}, and using Fubini's theorem to
   rewrite the second term of \eqref{eq: V}, we get
   \begin{align*}
   V^k(x,y) & = \frac{k}{\theta} \int_{[0,x]} F_0(t,y) dt
           - { \int \int}_{0 \le s \le t \le x} F_0(s,y) \frac{k (k-1)}{\theta^2}
           \left (1 - \frac{t-s}{\theta} \right )^{k-2} ds dt \\
   & =  \frac{k}{\theta} \int_{[0,x]} F_0(s,y) \left(
                1 - \frac{x-s}{\theta} \right )^{k-1} ds
    =  \int_{[0,x]} F_0(s,y) d Q_x^k (s),
   \end{align*}
   where, for $s\le x$,
   \begin{align*}
     Q_x^k (s) & = \int_0^s \frac{k}{\theta} \left(1-
     \frac{x-r}{\theta}\right)^{k-1} dr =   \left ( 1 - \frac{x-s}{\theta} \right )^k - \left (1 -
     \frac{x}{\theta} \right )^k.
   \end{align*}
   Thus, as $k\to \infty$, $Q_x^k(s)$ converges weakly to the
   distribution function with mass $1$ at $x$. Plugging in $y=\infty$
   in $V^{k}(x,y)$ yields $V_X^k(x) = \int_{[0,x]} F_{0X}(s)
   dQ_x^k(s)$. Furthermore, plugging in the expressions for $V_X^k$ and
   $G_k$ in \eqref{eq: H appendix} gives
   \begin{align*}
    H^k(x) & =  \int_{[0,x]} F_{0X}(s) dQ_x^k (s)
                +  \int_{[0,x]} (1-F_{0X}(s))  \frac{k}{\theta} \left(\frac{s}{\theta}\right)^{k-1}
                ds.
f   \end{align*}
   Hence, for $x<\theta$ we have $V^k (x,y) \rightarrow F_0(x,y)$,
   $V_X^k (x) \rightarrow F_{0X} (x)$ and $1-H^k (x) \rightarrow
   1-F_{0X}(x)$ as $k \rightarrow \infty$ for continuity points of the
   limits. The corollary then follows from the extended Helly-Bray
   theorem. \qed
\end{proof_of}

\newpage

\end{document}